\documentclass[12pt]{article}

\usepackage{latexsym,amsfonts,amsmath,amssymb,amsthm}
\usepackage[english]{babel}
\usepackage[latin1]{inputenc}
\usepackage[T1]{fontenc}
\usepackage{url}
\usepackage{graphicx}

\newtheorem{lemm}{Lemma}
\newtheorem{prop}{Proposition}
\newtheorem{theo}{Theorem}
\newtheorem*{ack}{Acknowledgments}

\newcommand{\R}{\ensuremath{\mathbb{R}}}
\renewcommand{\P}{\ensuremath{\mathbf{P}}}
\newcommand{\E}{\ensuremath{\mathbf{E}}}
\newcommand{\defeg}{\ensuremath{\overset{\hbox{\tiny{def}}}{=}}}
\newcommand{\V}{\ensuremath{\mathbb{V}}}
\newcommand{\A}{\ensuremath{\mathbb{A}}}
\newcommand{\T}{\ensuremath{\mathbb{T}}}
\newcommand{\inv}{\hbox{\tiny{-1}}}
\renewcommand{\S}{\ensuremath{\mathbb{S}}}

\begin{document}

\title{\textbf{A slow transient diffusion in a drifted stable potential}}

\author{\normalsize{\textsc{Arvind Singh
\footnote{\it{Laboratoire de Probabilit\'es
et Mod\`eles Al\'eatoires, Universit\'e Pierre et Marie Curie, 175 rue du
Chevaleret, 75013 Paris, France.} }}}}
\date{}
\maketitle \vspace{-1cm}
\begin{center}
Universit\'e Paris VI
\end{center}

\vspace*{0.2cm}

\begin{abstract}
We consider a diffusion process $X$ in a random potential $\V$ of
the form $\V_x = \S_x -\delta x$ where $\delta$ is a positive drift
and $\S$ is a strictly stable process of index $\alpha\in (1,2)$
with positive jumps. Then the diffusion is transient and $X_t /
\log^\alpha t$ converges in law towards an exponential distribution.
This behaviour contrasts with the case where $\V$ is a drifted
Brownian motion and provides an example of a transient diffusion in
a random potential which is as "slow" as in the recurrent setting.
\end{abstract}

\bigskip
{\small{


\bigskip
\noindent{\bf e-mail. }arvind.singh@ens.fr }}

\section{Introduction}
Let $(\V(x),\; x\in\R)$ be a two-sided stochastic process defined on
some probability space $(\Omega,\mathcal{F},\P)$. We call a
diffusion in the random potential $\V$ an informal solution $X$ of
the S.D.E:
\begin{equation*}
\left\{
\begin{array}{l}
dX_t = d\beta_t -\frac{1}{2}\V'(X_t)dt\\
X_0 = 0,
\end{array}
\right.
\end{equation*}
where $\beta$ is a standard Brownian motion independent of $\V$. Of
course, the process $\V$ may not be differentiable (for example when
$\V$ is a Brownian motion) and we should formally consider $X$ as a
diffusion whose conditional generator given $\V$ is
\begin{equation*}
\frac{1}{2}e^{\V(x)}\frac{d}{dx}\left(e^{-\V(x)}\frac{d}{dx}\right).
\end{equation*}
Such a diffusion may be explicitly constructed from a Brownian
motion through a random change of time and a random change of scale.
This class of processes has been widely studied for the last twenty
years and bears a close connection with the model of the random walk
in random environment (RWRE), see \cite{Zeitouni1} and \cite{Shi1}
for a survey on RWRE and \cite{Schumacher1}, \cite{Shi1} for the
connection between the two models.

This model exhibits many interesting features. For instance, when
the potential process $\V$ is a Brownian motion, the diffusion $X$
is recurrent and Brox \cite{Brox1} proved that $X_t/\log^2 t$
converges to a non-degenerate distribution. Thus, the diffusion is
much "slower" than in the trivial case $\V = 0$ (then $X$ is simply
a Brownian motion).

We point out that Brox's theorem is the analogue of Sinai's famous
theorem for RWRE \cite{Sinai1} (see also \cite{Golosov1} and
\cite{Kesten2}). Just as for the RWRE, this result is a consequence
of a so-called "localization phenomena": the diffusion is trapped in
some valleys of its potential $\V$. Brox's theorem may also be
extended to a wider class of potentials. For instance, when $\V$ is
a strictly stable process of index $\alpha\in(0,2]$, Schumacher
\cite{Schumacher1} proved that
\begin{equation*}
\frac{X_t}{\log^\alpha
t}\underset{t\to\infty}{\overset{\hbox{\tiny{law}}}{\longrightarrow}}b_\infty,
\end{equation*}
where $b_\infty$ is a non-degenerate random variable, whose
distribution depends on the parameters of the stable process $\V$.

There is also much interest concerning the behaviour of $X$ in the
transient case. When the potential is a drifted Brownian motion
\emph{i.e.} $\V_x = \mathbb{B}_x -\frac{\kappa}{2}x$ where
$\mathbb{B}$ is a two-sided Brownian motion and $\kappa>0$, then the
associated diffusion $X$ is transient toward $+\infty$ and its rate
of growth is polynomial and depends on $\kappa$. Precisely, Kawazu
and Tanaka \cite{KawazuTanaka1} proved that
\begin{itemize}
\item If $0 < \kappa <1$, then $\frac{1}{t^\kappa}X_t$ converges in law
towards a Mittag-Leffler distribution of index $\kappa$.
\item If $\kappa =1$, then $\frac{\log t}{t}X_t$ converges in
probability towards $\frac{1}{4}$.
\item If $\kappa >1$, then $\frac{1}{t}X_t$ converges almost surely towards $\frac{\kappa - 1}{4}$.
\end{itemize}
In particular, when $\kappa <1$, the rate of growth of $X$ is
sub-linear. Refined results on the rates of convergence for this
process were later obtained by Tanaka \cite{Tanaka2} and Hu \emph{et al.} \cite{HuShiYor1}.

In fact, this behaviour is not specific to diffusions in a drifted
Brownian potential. More generally, it is proved in \cite{Singh2}
that if $\V$ is a two-sided L\'evy process with no positive jumps
and if there exists $\kappa>0$ such $\E[e^{\kappa\V_1}] = 1$, then
the rate of growth of $X_t$ is linear when $\kappa > 1$ and of order
$t^\kappa$ when $0<\kappa < 1$ (see also \cite{Carmona1} for a law
of large numbers in a general L\'evy potential). These results are
the analogues of those previously obtained by Kesten \emph{et al.}
\cite{KestenKozlovSpitzer1} for the discrete model of the RWRE.

In this paper, we study the asymptotic behaviour of a diffusion in a
drifted stable potential. Precisely, let $(\S_x,\; x\in \R)$ denote
a two-sided c\`adl\`ag stable process with index $\alpha\in(1,2)$.
By two-sided, we mean that
\begin{enumerate}
\item[(a)] The process $(\S_x,\; x \geq 0)$ is strictly stable with
index $\alpha\in (1,2)$, in particular $\S_0=0$.
\item[(b)] For all $x_0\in\R$, the process $(\S_{x+x_0} - \S_{x_0},\; x\in \R)$ has the
same law as $\S$.
\end{enumerate}
It is well known that the L\'evy measure $\Pi$ of $\S$ has the form
\begin{equation}\label{RefLevyMeasure}
\Pi(dx) = \left(c^{+}\mathbf{1}_{\{x > 0\}} +c^{-}\mathbf{1}_{\{x <
0\}}\right)\frac{dx}{|x|^{\alpha+1}}
\end{equation}
where $c^+$ and $c^-$ are two non-negative constants such that $c^+
+ c^{-} >0$. In particular, the process $(\S_{x}\, , \, x\geq 0)$
has no positive jumps (resp. no negative jumps) if and only if $c^+
=0$ (resp. $c^- = 0$). Given $\delta >0$, we consider a diffusion
$X$ in the random potential
\begin{equation*}
\V_x =\S_x - \delta x.
\end{equation*}
Since the index $\alpha$ of the stable process $\S$ is larger than
$1$, we have $\E[\V_x] = -\delta x$, and therefore
\begin{equation*}
\lim_{x\to+\infty}\V_x = -\infty \quad\hbox{and}\quad \lim_{x\to-\infty}\V_x = +\infty\quad\hbox{almost surely.}
\end{equation*}
This implies that the random diffusion $X$ is transient toward
$+\infty$. We already mentioned that, when $\S$ has no positive
jumps (\emph{i.e.} $c^+=0$), the rate of transience of $X$ is given
in \cite{Singh2} and $X_t$ has a polynomial growth. Thus, we now
assume that $\S$ possesses positive jumps.
\begin{theo} \label{MainTheo} Assume that $c^+ >0$, then
\begin{equation*}
\frac{X_t}{\log^\alpha
t}\underset{t\to\infty}{\overset{\hbox{\textup{\tiny{law}}}}{\longrightarrow}}
\mathcal{E}\left(\frac{c^+}{\alpha}\right),
\end{equation*}
where $\mathcal{E}(c^+/\alpha)$ denotes an exponential law with
parameter $c^+/\alpha$. This result also holds with $\sup_{s\leq
t}X_s$ or $\inf_{s\geq t}X_s$ in place of $X_t$.
\end{theo}
The asymptotic behaviour of $X$ is in this case very different from
the one observed when $\V$ is a drifted Brownian motion. Here, the
rate of growth is very slow: it is the same as in the recurrent
setting. We also note that neither the rate of growth nor the
limiting law depend on the value of the drift parameter $\delta$.

Theorem \ref{MainTheo} has a simple heuristic explanation:  the
"localisation phenomena" for the diffusion $X$ tells us that the
time needed to reach a positive level $x$ is approximatively
exponentially proportional to the biggest ascending barrier of $\V$
on the interval $[0,x]$. In the case of a Brownian potential, or
more generally a spectrally negative L\'evy potential, the addition
of a negative drift somehow "kills" the ascending barriers, thus
accelerating the diffusion and leading to a polynomial rate of
transience. However, in our setting, the biggest ascending barrier
on $[0,x]$ of the stable process $\S$ is of the same order as its
biggest jump on this interval. Since, the addition of a drift has no
influence on the jumps of the potential process, the time needed to
reach level $x$ still remains of the same order as in the recurrent
case (\emph{i.e.} when the drift is zero) and yields a logarithmic
rate of transience.

\section{Proof of the theorem.}

\subsection{Representation of $X$ and of its hitting times.}

In the remainder of this paper, we indifferently use the notation $\V_x$ or $\V(x)$.
Let us first recall the classical representation of the diffusion
$X$ in the random potential $\V$ from a Brownian motion through a
random change of scale and a random change of time (see \cite{Brox1}
or \cite{Shi1} for details). Let $(B_t,\; t\geq 0)$ denote a
standard Brownian motion independent of $\V$ and let $\sigma$ stand
for its hitting times:
\begin{equation*}
\sigma(x) \defeg \inf(t\geq 0\, , \, B_t = x).
\end{equation*}
Define the scale function of the diffusion $X$,
\begin{equation}\label{defA}
\A(x) \defeg \int_{0}^{x}e^{\V_y}dy\quad\hbox{ for $x\in\R$.}
\end{equation}
Since $\lim_{x\to +\infty}\V_x/x = -\delta$ and $\lim_{x\to
-\infty}\V_x/x = \delta$ almost surely, it is clear that
\begin{equation*}
\A(\infty) = \lim_{x\to +\infty}\A(x)
<\infty\quad \hbox{ and }\quad \lim_{x\to-\infty}\A(x) = -\infty\quad\hbox{ almost surely.}
\end{equation*}
Let $\A^{\inv} : (-\infty,\A(\infty))\mapsto\R$ denote the inverse of
$\A$ and define
\begin{equation*}
\T(t) \defeg \int_{0}^{t}e^{-2\V(\A^{\inv}(B_s))}ds\quad\hbox{ for $0\leq
t < \sigma(\A(\infty))$.}
\end{equation*}
Similarly, let $\T^{\inv}$ denote the inverse of $\T$. According to
Brox \cite{Brox1} (see also \cite{Shi1}), the diffusion $X$ in the
random potential $\V$ may be represented in the form
\begin{equation}\label{representationX}
X_t = \A^{\inv}\Big(B_{\T^{\inv}(t)}\Big).
\end{equation}
It is now clear that, under our assumptions, the diffusion $X$ is
transient toward $+\infty$. We will study $X$ via its
hitting times $H$ defined by
\begin{equation*}
H(r) \defeg \inf(t\geq 0\, , \, X_t = r) \quad \hbox{ for $r\geq0$.}
\end{equation*}
Let $(L(t,x),\; t\geq 0,\; x\in\R)$ stand for the bi-continuous
version of the local time process of $B$. In view of
(\ref{representationX}), we can write
\begin{equation*}
H(r) = \T\left(\sigma(\A(r))\right) =
\int_{0}^{\sigma(\A(r))}e^{-2\V(\A^{\inv}(B_s))}ds =
\int_{-\infty}^{\A(r)}e^{-2\V(\A^{\inv}(x))}L(\sigma(\A(r)),x)dx.
\end{equation*}
Making use of the change of variable $x = \A(y)$, we get
\begin{equation}\label{defH}
H(r) = \int_{-\infty}^{r}e^{-\V_y}L(\sigma(\A(r)),\A(y))dy = I_1(r) + I_2(r)
\end{equation}
where
\begin{eqnarray*}
I_1(r) &\defeg& \int_{0}^{r}e^{-\V_y}L(\sigma(\A(r)),\A(y))dy,\\
I_2(r) &\defeg& \int_{-\infty}^{0}e^{-\V_y}L(\sigma(\A(r)),\A(y))dy.\\
\end{eqnarray*}

\subsection{Proof of Theorem \ref{MainTheo}.}
Given a c\`adl\`ag process $(Z_t\hbox{ , }t\geq 0)$, we denote by
$\Delta_t Z = Z_t - Z_{t-}$ the size of the jump at time $t$. We
also use the notation $Z^\natural_t$ to denote the largest positive
jump of $Z$ before time $t$,
\begin{equation*}
Z^\natural_t \defeg \sup(\Delta_s,\;0\leq s\leq t).
\end{equation*}
Let $Z^\#_t$ stand for the largest ascending barrier before time $t$, namely:
\begin{equation*}
Z^\#_t \defeg \sup_{0\leq x \leq y \leq t}(Z_y - Z_x).
\end{equation*}
We also define the functionals:
\begin{equation*}
\begin{aligned}
\overline{Z}_t &\defeg \sup_{s\in[0,t]}Z_s&\quad\hbox{(running unilateral maximum)}\\
\underline{Z}_t &\defeg \inf_{s\in[0,t]}Z_s&\quad\hbox{(running unilateral minimum)}\\
Z^*_t&\defeg\sup_{s\in[0,t]}|Z_s|&\quad\hbox{(running bilateral
supremum)}
\end{aligned}
\end{equation*}
We start with a simple lemma concerning the fluctuations of the
potential process.
\begin{lemm}\label{lemmfluc}
There exist two constants $c_1,c_2>0$ such that for all $a,x>0$
\begin{equation}\label{flucleq}
\P\{\V^\#_x \leq a\}\leq e^{-c_1\frac{x}{a^\alpha}},
\end{equation}
and whenever $\frac{a}{x}$ is sufficiently large,
\begin{equation}\label{flucgeq}
\P\{\V^*_x > a\}\leq c_2\frac{x}{a^\alpha}.
\end{equation}
\end{lemm}
\begin{proof}
Recall that $\V_x = \S_x - \delta x$. In view of the form of the
density of the L\'evy measure of $\S$ given in (\ref{RefLevyMeasure}),
we get
\begin{equation*}
\P\{\V^\#_x \leq a\}\leq \P\{\V^\natural_x \leq a\} =
\exp\left(-x\int_{a}^{\infty}\frac{c^+}{y^\alpha}dy\right) =
\exp\left(-\frac{c^+}{\alpha}\frac{x}{a^\alpha}\right).
\end{equation*}
This yields (\ref{flucleq}). From the scaling property of the stable
process $\S$, we also have
\begin{equation*}
\P\{\V^*_x > a\} =
\P\left\{x^{\frac{1}{\alpha}}\sup_{t\in[0,1]}|\S_t - \delta
x^{1-\frac{1}{\alpha}}t|> a\right\} \leq \P\left\{\S^*_1>
\frac{a}{x^{\frac{1}{\alpha}}} - \delta
x^{1-\frac{1}{\alpha}}\right\}.
\end{equation*}
Notice further that $a/x^{1/\alpha} - \delta x^{1-1/\alpha} >
a/(2x^{1/\alpha})$ whenever $a/x$ is large enough, therefore, making
use of a classical estimate concerning the tail distribution of the
stable process $\S$ (\emph{c.f.} Proposition $4$, p221 of
\cite{Bertoin1}), we find that
\begin{equation*}
\P\{\V^*_x > a\}\leq \P\left\{\S^*_1>
\frac{a}{2x^{\frac{1}{\alpha}}}\right\}\leq \P\left\{\overline{\S}_1
> \frac{a}{2x^{\frac{1}{\alpha}}}\right\}
+\P\left\{\underline{\S}_1> \frac{a}{2x^{\frac{1}{\alpha}}}\right\}
\leq c_2\frac{x}{a^\alpha}.
\end{equation*}

\end{proof}
\begin{prop}\label{lem1}There exists a constant $c_3>0$ such that,
for all $r$ sufficiently large and all $x\geq 0$,
\begin{equation*}
\P\{\V^\#_{r} \geq x + \log^4 r\} - c_3 e^{-\log^2 r} \leq
\P\{\log I_1(r) \geq x\} \leq \P\{\V^\#_r \geq x - \log^4 r\} + c_3
e^{-\log^2 r}.
\end{equation*}
\end{prop}
\begin{proof}
This estimate was first proved by Hu and Shi (see Lemma 4.1 of
\cite{HuShi1}) when the potential process is close to a standard
Brownian motion. A similar result is given in Proposition 3.2 of
\cite{Singh1} when $\V$ is a random walk in the domain of attraction
of a stable law. As explained by Shi \cite{Shi1}, the key idea is
the combined use of Ray-Knight's Theorem and Laplace's method.
However, in our setting, additional difficulties appear since the
potential process is neither flat on integer interval nor
continuous. We shall therefore give a complete proof but one can
still look in \cite{HuShi1} and \cite{Singh1} for additional
details. Recall that
\begin{equation*}
I_1(r) = \int_{0}^{r}e^{-\V_y}L(\sigma(\A(r)),\A(y))dy,\\
\end{equation*}
where $L$ is the local time of the Brownian motion $B$ (independent
of $\V$). Let $(U(t),\; t\geq 0)$ denote a two-dimensional squared
Bessel process starting from zero, also independent of $\V$.
According to the first Ray-Knight Theorem (\emph{c.f.} Theorem 2.2
p455 of \cite{RevuzYor1}), for any $x>0$ the process
$(L(\sigma(x),x-y),\; 0\leq y\leq x)$ has the same law as
$(U(y),\;0\leq y\leq x)$. Therefore, making use of the scaling
property of the Brownian motion and the independence of $\V$ and
$B$, for each fixed $r>0$, the random variable $I_1(r)$ has the same
law as
\begin{equation*}
\widetilde{I}_1(r) \defeg \A(r)\int_{0}^{r}e^{-\V_y}U\left(\frac{\A(r)-\A(y)}{\A(r)}\right)dy.
\end{equation*}
We simply need to prove the proposition for $\widetilde{I}_1$
instead of $I_1$. In the rest of the proof, we assume that $r$ is
very large. We start with the upper bound. Define the event
\begin{equation*}
\mathcal{E}_1\defeg \left\{\sup_{t\in(0,1]}\frac{U(t)}{t \log\left(\frac{8}{t}\right)}\leq r\right\}.
\end{equation*}
According to Lemma 6.1 of \cite{HuShi1}, $\P\{\mathcal{E}_1^c\}\leq c_4
e^{-r/2}$ for some constant $c_4>0$. On $\mathcal{E}_1$, we have
\begin{eqnarray*}
\widetilde{I}_1(r) &\leq& r\int_{0}^{r}e^{-\V_y}(\A(r)-\A(y))\log\left(\frac{\A(r)}{\A(r)-\A(y)}\right)dy\\
&=&r\int_{0}^{r}\left(\int_{y}^{r}e^{\V_z-\V_y}dz\right)\log\left(\frac{\A(r)}{\A(r)-\A(y)}\right)dy\\
&\leq&r^{2}e^{\V^\#_r}\int_{0}^{r}\log\left(\frac{\A(r)}{\A(r)-\A(y)}\right)dy.
\end{eqnarray*}
Notice also that $\A(r) = \int_{0}^{r}e^{\V_z}dz\leq
re^{\overline{\V}_r}$ and similarly $\A(r)-\A(y) \geq
(r-y)e^{\underline{\V}_r}$. Therefore
\begin{eqnarray*}
\int_{0}^{r}\log\left(\frac{\A(r)}{\A(r)-\A(y)}\right)dy
&\leq& r(\overline{\V}_r - \underline{\V}_r) + \int_{0}^{r}\log\left(\frac{8r}{r-y}\right)dy\\
&=& r(\overline{\V}_r - \underline{\V}_r + 1 + \log 8).
\end{eqnarray*}
Define the set $\mathcal{E}_2 \defeg \{\overline{\V}_r - \underline{\V}_r
\leq e^{\log^3 r}\}$. In view of Lemma \ref{lemmfluc},
\begin{equation*}
\P\{ \mathcal{E}_2^c\} \leq \P\left\{ \V^*_r > \frac{1}{2}e^{\log^3
r}\right\} \leq e^{-\log^2 r}.
\end{equation*}
Therefore, $\P\{(\mathcal{E}_1\cap\mathcal{E}_2)^c\}\leq 2e^{-\log^2
r}$ and on $\mathcal{E}_1\cap\mathcal{E}_2$,
\begin{equation*}
\widetilde{I}_1(r) \leq r^3 (e^{\log^3 r} + 1 + \log 8)e^{\V^\#_r}
\leq e^{\log^4 r + \V^\#_r}.
\end{equation*}
This completes the proof of the upper bound. We now deal with the
lower bound. Define the sequence $(\gamma_k,\;k\geq 0)$ by induction
\begin{equation*}
\left\{
\begin{array}{cll}
\gamma_0 &\defeg& 0, \\
\gamma_{k+1} &\defeg& \inf(t>\gamma_n,\; |\V_t - \V_{\gamma_k}|\geq 1).
\end{array}\right.
\end{equation*}
The sequence $(\gamma_{k+1}-\gamma_k,\; k\geq 0)$ is \emph{i.i.d.} and
distributed as $\gamma_1 = \inf(t> 0,\; |\V_t|\geq 1)$. We denote by
$[x]$ the integer part of $x$. We also use the notation $\epsilon \defeg
e^{-\log^3 r}$. Consider the following events
\begin{eqnarray*}
\mathcal{E}_3 &\defeg& \left\{ \gamma_{[r^2]} > r\right\},\\
\mathcal{E}_4 &\defeg& \left\{ \gamma_{k} - \gamma_{k-1} \geq
2\epsilon\;\hbox{ for all $k=1,2\ldots,[r^2]$}\right\}.
\end{eqnarray*}
In view of Cramer's large deviation Theorem and since $r$ is very
large, we get that $\P\{\mathcal{E}_3^c\}\leq e^{-r}$. We also have
\begin{eqnarray*}
\P\{\mathcal{E}_4^c \} \leq \sum_{k=1}^{[r^2]} \P\{ \gamma_{k} -
\gamma_{k-1} < 2\epsilon\} &\leq& [r^2]\P\{\gamma_1 < 2\epsilon\}\\
&\leq&[r^2]\P\{\V^*_{2\epsilon} \geq 1\}\\
&\leq& e^{-\log^2 r},
\end{eqnarray*}
where we used Lemma \ref{lemmfluc} for the last inequality. Define
also
\begin{equation*}
\mathcal{E}_{5} \defeg \{ |\V_x - \V_r|< 1\;\hbox{ for all
$x\in[r-2\epsilon,r]$}\}.
\end{equation*}
From time reversal, the processes $(\V_{t},\;0\leq t\leq 2\epsilon)$
and $(\V_{r} -\V_{(r-t)^-},0\leq t\leq 2\epsilon)$ have the same
law. Thus,
\begin{equation*}
\P\{\mathcal{E}_5^c\}\leq \P\{\V^*_{2\epsilon} \geq 1\} \leq
e^{-\log^2 r}.
\end{equation*}
Setting $\mathcal{E}_6 \defeg
\mathcal{E}_3\cap\mathcal{E}_4\cap\mathcal{E}_5$, we get
$\P\{\mathcal{E}_6^c\}\leq 3e^{-\log^2 r}$. Moreover, it is easy to
check (see figure $1$) that on $\mathcal{E}_6$, we can always find
$x_-,x_+$ such that:
\begin{equation*}
\left\{
\begin{array}{l}
0 \leq x_- \leq x_+ \leq r - 2\epsilon, \\
\hbox{for any $a\in[x_-,x_- + \epsilon]$, $|\V_{x_-}-\V_a| \leq 2$,}\\
\hbox{for any $b\in[x_+,x_+ + \epsilon]$, $|\V_{x_+}-\V_b| \leq 2$,}\\
\V_{x_+} - \V_{x_-} \geq \V^\#_r - 4.
\end{array}\right.
\end{equation*}
\begin{figure}[top]
\setlength\unitlength{1cm}
\begin{picture}(15,10)
\put(1,0){\includegraphics[height=10cm,angle=0]{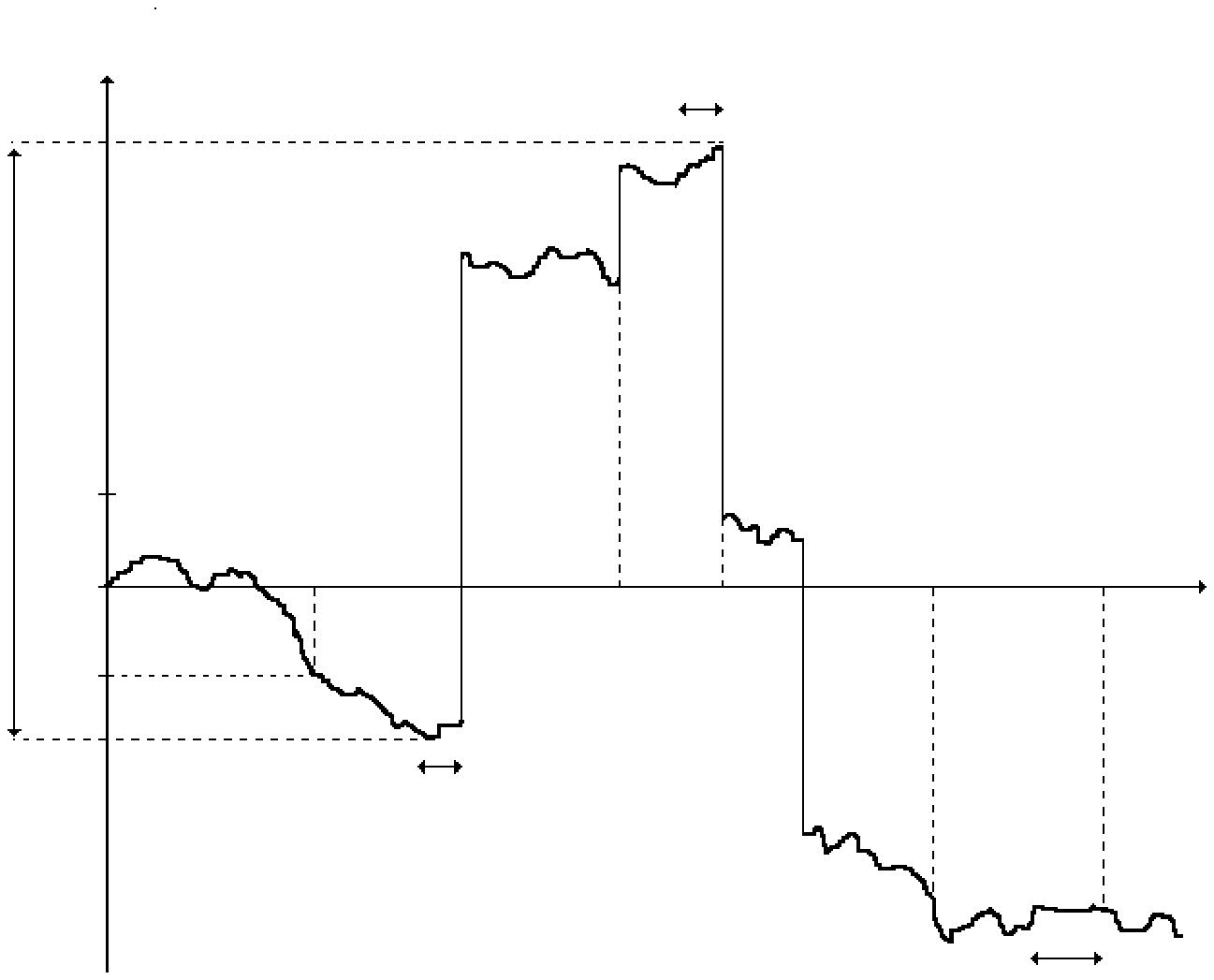}}
\put(2.1,4.4){$0$} \put(2.1,5.3){$1$} \put(1.8,3.6){$-1$}
\put(0.9,5.8){$\V^\#_r$} \put(4.2,4.7){$\gamma_1$}
\put(5.6,4.7){$\gamma_2$}\put(10,4.7){$\gamma_6$}
\put(7.1,4.2){$\gamma_3$} \put(8,4.2){$\gamma_4$}
\put(8.8,4.2){$\gamma_5$} \put(11.7,4.65){$r$} \put(4.95,2.5){$x_-$}
\put(5.85,2.5){$x_{-}\!\!+\!\epsilon$} \put(7.35,9.2){$x_+$}
\put(8.25,9.2){$x_{+}\!\!+\!\epsilon$}
\put(10.4,0.7){$r\!\!-\!\!2\epsilon$} \put(11.8,0.7){$r$}
\end{picture} \caption{Sample path of $\V$ on $\mathcal{E}_6$.}
\end{figure}
Let us also define
\begin{eqnarray*}
\mathcal{E}_7 &\defeg& \mathcal{E}_6\cap\left\{
\inf_{y\in[x_-,x_-+\epsilon]}U\left(\frac{\A(r)-\A(y)}{\A(r)}\right)\geq
\frac{\A(r)-\A(x_-)}{\A(r)}e^{-2\log^2 r}\right\},\\
\mathcal{E}_8 &\defeg& \left\{\V^\natural_{r}\geq 3\log^2 r\right\}.
\end{eqnarray*}
We finally set $\mathcal{E}_9 \defeg \mathcal{E}_7\cap\mathcal{E}_8$.
Then on $\mathcal{E}_9$, we have, for all $r$ large enough,
\begin{eqnarray*}
\widetilde{I}_1(r) &\geq&
\A(r)\int_{x_-}^{x_-+\epsilon}e^{-\V_y}U\left(\frac{\A(r)-\A(y)}{\A(r)}\right)dy\\
&\geq& e^{-\V_{x_-}-2 -2\log^2 r}\int_{x_-}^{x_- +\epsilon}(\A(r) -
\A(x_-))dy\\
&=&e^{-\V_{x_-}-2 -2\log^2 r -\log^3 r}\int_{x_-}^r e^{\V_y}dy\\
&\geq&e^{-\V_{x_-}-2 -2\log^2 r -\log^3 r}\int_{x_+}^{x_+ +
\epsilon}
e^{\V_{y}}dy\\
&\geq&e^{\V_{x_+}-\V_{x_-}-4 -2\log^2 r - 2\log^3 r}\\
&\geq&e^{\V^\#_r - \log^4 r}.
\end{eqnarray*}
This proves the lower bound on $\mathcal{E}_9$. It simply remains to
show that $\P\{\mathcal{E}_9^c\}\leq c_5 e^{-\log^2 r}$. According
to Lemma $6.1$ of \cite{HuShi1}, for any $0<a<b$ and any $\eta>0$,
we have
\begin{equation*}
\P\left\{ \inf_{a<t<b}U(t) \leq \eta b\right\}\leq 2\sqrt{\eta} +
2\exp\left(-\frac{\eta}{2(1-a/b)}\right).
\end{equation*}
Therefore, making use of the independence of $\V$ and $U$, we find
\begin{eqnarray*}
\P\{\mathcal{E}_9^c\}&\leq& \P\{\mathcal{E}_6^c\} +
\P\{\mathcal{E}_8^c\} + \P\{\mathcal{E}_7^c \cap \mathcal{E}_6 \cap \mathcal{E}_8\} \\
&\leq& \P\{\mathcal{E}_6^c\} +
\P\{\mathcal{E}_8^c\} + 2e^{-\log^2 r}+2\E\left[e^{\frac{1}{2}\mathbb{J}(r)e^{-2\log^2
r}}\mathbf{1}_{\mathcal{E}_6 \cap \mathcal{E}_8}\right],
\end{eqnarray*}
where
\begin{equation*}
\mathbb{J}(r) \defeg \frac{\A(r) - \A(x_-)}{\A(x_-+\epsilon)-\A(x_-)}.
\end{equation*}
We have already proved that $\P\{\mathcal{E}_6^c\}\leq 3e^{-\log^2
r}$. Using Lemma \ref{lemmfluc}, we also check that
$\P\{\mathcal{E}_8^c\}\leq e^{-\log^2 r}$. Thus, it remains to show
that
\begin{equation}\label{RefInegE8}
\E\left[e^{\frac{1}{2}\mathbb{J}(r)e^{-2\log^2
r}}\mathbf{1}_{\mathcal{E}_6 \cap \mathcal{E}_8}\right] \leq c_6
e^{-\log^2 r}.
\end{equation}
Notice that, on $\mathcal{E}_6$,
\begin{equation*}
\A(r) - \A(x_-) = \int_{x_-}^{r}e^{\V_y}dy\geq \int_{x_+}^{x_+
+\epsilon}e^{\V_y}dy\geq e^{\log^3 r + \V_{x_+} - 2},
\end{equation*}
and also
\begin{equation*}
\A(x_- +\epsilon) - \A(x_-) =\int_{x_-}^{x_- +\epsilon} e^{\V_y}dy
\leq e^{\log^3 r + \V_{x_-} + 2}.
\end{equation*}
Therefore, on $\mathcal{E}_6\cap\mathcal{E}_8$,
\begin{equation*}
\mathbb{J}(r) \geq e^{\V_{x_+} - \V_{x_-} - 4} \geq e^{\V^\#_r - 8} \geq
e^{\V^\natural_r - 8}\geq  e^{3\log^2 r - 8}
\end{equation*}
which clearly yields (\ref{RefInegE8}) and the proof of the proposition
is complete.
\end{proof}

\begin{lemm}\label{lem2} We have
\begin{equation*}
\frac{\V^\#_r}{r^{1/\alpha}}\underset{r\to\infty}{\overset{\hbox{\tiny{law}}}{\longrightarrow}}\S^\natural_1.
\end{equation*}
\end{lemm}

\begin{proof}
Let $f:[0,1]\mapsto\R$ be a deterministic c\`adl\`ag function. For
$\lambda\geq 0$, define
\begin{equation*}
f_\lambda (x) \defeg f(x) - \lambda x.
\end{equation*}
We first show that
\begin{equation}\label{RefCVf}
\lim_{\lambda\to\infty}f^{\#}_{\lambda}(1)  = f^{\natural}(1).
\end{equation}
It is clear that $f^{\natural}(1) = f_{\lambda}^{\natural}(1) \leq
f_{\lambda}^{\#}(1)$ for any $\lambda>0$. Thus, we simply need to
prove that $\limsup f^{\#}_\lambda(1) \leq f^{\natural}(1)$. Let
$\eta>0$ and set
\begin{eqnarray*}
A(\eta,\lambda) &\defeg& \sup\left(f_{\lambda}(y) - f_{\lambda}(x)\hbox{ , }0\leq x\leq y\leq 1\hbox{ and }y-x \leq \eta\right),\\
B(\eta,\lambda) &\defeg& \sup\left(f_{\lambda}(y) - f_{\lambda}(x)\hbox{ , }0\leq x\leq y\leq 1\hbox{ and }y-x > \eta\right),
\end{eqnarray*}
so that
\begin{equation}\label{refap1}
f^{\#}_{\lambda}(1) = \max(A(\eta,\lambda),B(\eta,\lambda)).
\end{equation}
Notice that $A(\eta,\lambda) \leq A(\eta)$ where
\begin{equation*}
A(\eta) \defeg A(\eta,0) = \sup\left(f(y) - f(x)\hbox{ , }0\leq x\leq y\leq 1\hbox{ and }y-x \leq \eta\right).
\end{equation*}
Since $f$ is c\`adl\`ag, we have $\lim_{\eta\to 0} A(\eta) =
f^\natural(1)$. Thus, for any $\varepsilon>0$, we can find
$\eta_0>0$ small enough such that
\begin{equation}\label{refap2}
\limsup_{\lambda\to\infty} A(\eta_0,\lambda) \leq f^\natural(1)+\varepsilon.
\end{equation}
Notice also that
\begin{eqnarray*}
B(\eta_0,\lambda) & \leq & \sup\left(f(y) - f(x) - \eta_0\lambda\hbox{ , }0\leq x\leq y\leq 1\hbox{ and }y-x > \eta_0\right)\\
& \leq & f^\#(1) - \eta_0\lambda
\end{eqnarray*}
which implies
\begin{equation}\label{refap3}
\lim_{\lambda\to\infty}B(\eta_0,\lambda) =-\infty.
\end{equation}
The combination of (\ref{refap1}), (\ref{refap2}) and (\ref{refap3})
yield (\ref{RefCVf}). Making use of the scaling property of the
stable process $\S$, for any fixed $r>0$,
\begin{equation*}
(\V_y\hbox{ , } 0\leq y \leq r) \overset{\hbox{\tiny{law}}}{=}
(r^{1/\alpha}\S_y - \delta r y\hbox{ , }0\leq y\leq 1).
\end{equation*}
Therefore, setting $\mathbb{R}(z) = (\S_\cdot - z\cdot)^\#_1$, we
get the equality in law:
\begin{equation}\label{egallaw}
\frac{\V^\#_r}{r^{1/\alpha}}\overset{\hbox{\tiny{law}}}{=}\mathbb{R}(\delta r^{1-1/\alpha}).
\end{equation}
Making use of (\ref{RefCVf}), we see that $\mathbb{R}(z)$ converges
almost surely towards $\S^\natural_1$ as $z$ goes to infinity. Since
$\alpha>1$ and $\delta>0$, we also have $\delta r^{1-1/\alpha} \to
\infty$ as $r$ goes to infinity and we conclude from (\ref{egallaw})
that
\begin{equation*}
\frac{\V^\#_r}{r^{1/\alpha}}\underset{r\to\infty}{\overset{\hbox{\tiny{law}}}{\longrightarrow}}\S^\natural_1.
\end{equation*}
\end{proof}

\begin{proof}[Proof of Theorem \ref{MainTheo}.]
Recall that the random variable $\S^\natural_1$ denotes the largest
positive jump of $\S$ over the interval $[0,1]$. Thus, according to
the density of the L\'evy measure of $\S$,
\begin{equation}\label{refz1}
\P\{\S^\natural_1 \leq x\} =
\exp\left(-\int_{x}^{\infty}\frac{c^+}{y^{\alpha+1}}dy\right) =
\exp\left(-\frac{c^+}{\alpha y^\alpha}\right).
\end{equation}
On the one hand, the combination of Lemma \ref{lem1} and \ref{lem2}
readily shows that
\begin{equation}\label{refz2}
\frac{\log(I_1(r))}{r^{1/\alpha}}\underset{r\to\infty}{\overset{\hbox{\tiny{law}}}{\longrightarrow}}\S^\natural_1.
\end{equation}
On the other hand, the random variables $\A(\infty) =
\lim_{x\to\infty} \A(x)$ and $\int_{-\infty}^{0}e^{-\V_y}dy$ have
the same law. Moreover, we already noticed that these random
variables are almost surely finite. Since the function $L(t,\cdot)$
is, for any fixed $t$, continuous with compact support, we get
\begin{equation*}
I_2(r) = \int_{-\infty}^{0}e^{-\V_y}L(\sigma(\A(r)),\A(y))dy \leq
\sup_{z\in(-\infty,0]}L(\sigma(\A(\infty)),z)\int_{-\infty}^{0}e^{-\V_y}dy
< \infty.
\end{equation*}
Therefore,
\begin{equation}\label{refz3}
\sup_{r\geq 0}I_2(r) < \infty\quad\hbox{almost surely.}
\end{equation}
Combining (\ref{defH}), (\ref{refz2}) and (\ref{refz3}), we deduce that
\begin{equation*}
\frac{\log(H(r))}{r^{1/\alpha}}\underset{r\to\infty}{\overset{\hbox{\tiny{law}}}{\longrightarrow}}\S^\natural_1
\end{equation*}
which, from the definition of the hitting times $H$, yields
\begin{equation*}
\frac{\sup_{s\leq t}X_s}{\log^\alpha
t}\underset{t\to\infty}{\overset{\hbox{\tiny{law}}}{\longrightarrow}}\left(\frac{1}{\S^\natural_1}\right)^\alpha.
\end{equation*}
According to (\ref{refz1}), the random variable
$(1/\S^\natural_1)^\alpha$ has an exponential distribution with
parameter $c^+/\alpha$ so the proof of the theorem for $\sup_{s\leq
t}X_s$ is complete. We finally use the classical argument given by
Kawazu and Tanaka, p201 \cite{KawazuTanaka1} to obtain the
corresponding results for $X_t$ and $\inf_{s\geq t}X_s$.
\end{proof}

\begin{ack}
I would like to thank Yueyun Hu for his precious advices.
\end{ack}

\end{document}